\documentclass[11pt]{amsart}

\usepackage{amssymb,amscd,amsmath,amsthm}

\theoremstyle{plain}
\newtheorem{thm}{Theorem}[section]
\newtheorem{theorem}[thm]{Theorem}

\newtheorem{cor}[thm]{Corollary}

\theoremstyle{definition}

\newtheorem{prob}[thm]{Problem}

\theoremstyle{remark}
\newtheorem{rem}{Remark}

\newtheorem*{remark*}{Remark}

\newcommand{\cF}{{\mathcal{F}}}

\newcommand{\cH}{{\mathcal{H}}}

\newcommand{\cL}{{\mathcal{L}}}

        \newcommand{\field}[1]{{\mathbb{#1}}}
        
        \newcommand{\ZZ}{\field{Z}}
        \newcommand{\QQ}{\field{Q}}
        \newcommand{\RR}{\field{R}}

\newcommand{\im}{\operatorname{Im}}

\newcommand{\Tr}{\operatorname{Tr}}
\newcommand{\tr}{\operatorname{tr}}

\begin{document}

\title[Adiabatic limits and the spectrum of foliations] {Adiabatic limits and the spectrum of the
Laplacian on foliated manifolds}
\author{Yuri A. Kordyukov}
\address{Institute of Mathematics,
         Russian Academy of Sciences,
         112~Chernyshevsky str., 450077 Ufa, Russia} \email{yurikor@matem.anrb.ru}

\author{Andrey A. Yakovlev}
\address{Department of Mathematics, Ufa State Aviation Technical University, 12 K. Marx str.,
450000 Ufa, Russia} \email{yakovlevandrey@yandex.ru}
\thanks{Supported by the Russian Foundation of Basic Research
(grant no. 06-01-00208)}

\begin{abstract}
We present some recent results on the behavior of the spectrum of
the differential form Laplacian on a Riemannian foliated manifold
when the metric on the ambient manifold is blown up in directions
normal to the leaves (in the adiabatic limit).
\end{abstract}
\maketitle \tableofcontents

\section*{Introduction}
Let $(M,{\mathcal F})$ be a closed foliated manifold, $\dim M =
n$, $\dim {\mathcal F} = p$, $p+q=n$, endowed with a Riemannian
metric $g$. Then we have a decomposition of the tangent bundle to
$M$ into a direct sum $TM=F\oplus H$, where $F=T{\mathcal F}$ is
the tangent bundle to $\cF$ and $H=F^{\bot}$ is the orthogonal
complement of $F$, and the corresponding decomposition of the
metric: $ g=g_{F}+g_{H}$. Define a one-parameter family $g_{h}$ of
Riemannian metrics on $M$ by
\begin{equation}\label{e:gh}
g_{h}=g_{F} + {h}^{-2}g_{H}, \quad h > 0.
\end{equation}
By the adiabatic limit, we will mean the asymptotic behavior of
Riemannian manifolds $(M, g_h)$ as $h\to 0$.

In this form, the notion of the adiabatic limit was introduced by
Witten \cite{Witten85} in the study of the global anomaly. He
considered a family of Dirac operators acting along the fibers of
a Riemannian fiber bundle over the circle and gave an argument
relating the holonomy of the determinant line bundle of this
family to the adiabatic limit of the eta invariant of the Dirac
operator on the total space. Witten's result was proved rigorously
in \cite{BismutFreedI}, \cite{BismutFreedII} and \cite{Cheeger87},
and extended to general Riemannian bundles in \cite{BismutCheeger}
and \cite{Dai91}. This study gave rise to the development of
adiabatic limit technique for analyzing the behavior of certain
spectral invariants under degeneration that has many applications
in the local index theory (see, for instance, \cite{BismutICM98}).

New properties of adiabatic limits were discovered by Mazzeo and
Melrose \cite{MazzeoMelrose}. They showed that in the case of a
fibration, a Taylor series analysis of so called small eigenvalues
in the adiabatic limit and the corresponding eigenforms leads
directly to a spectral sequence, which is isomorphic to the Leray
spectral sequence. This result was used in \cite{Dai91}, and
further developed in \cite{Forman95}, where the very general
setting of any pair of complementary distributions is considered.
Nevertheless, the most interesting results of \cite{Forman95} are
only proved for foliations satisfying very restrictive conditions.
The ideas from \cite{MazzeoMelrose} and \cite{Forman95} were also
applied in the case of the contact-adiabatic (or sub-Riemannian)
limit in \cite{Ge95,Rumin00}.

In this paper, we will discuss extensions of the results mentioned
above to the adiabatic limits on foliated manifolds. For any
$h>0$, we will consider the Laplace operator $\Delta_{h}$ on
differential forms defined by the metric $g_h$. It is a
self-adjoint, elliptic, differential operator with the positive,
scalar principal symbol in the Hilbert space $L^2(M,\Lambda
T^{*}M,g_h)$ of square integrable differential forms on $M$,
endowed with the inner product induced by $g_h$, which has
discrete spectrum. Denote by
\[
0\leq\lambda_0(h)\leq\lambda_1(h)\leq\lambda_2(h)\leq\cdots
\]
the spectrum of $\Delta_{h}$, taking multiplicities into account.

In Section~\ref{s:adiab}, we discuss the asymptotic behavior as
$h\to 0$ of the trace of $f(\Delta_h)$:
\[
\tr f(\Delta_h)=\sum_{i=0}^{+\infty}f(\lambda_i(h)),
\]
for any sufficiently nice function $f$, say, for $f\in S(\RR)$.
The results given in this Section should be viewed as a very first
step in extending the adiabatic limit technique to analyze the
behavior of spectral invariants to the case of foliations.

In Section~\ref{s:spec} we study ``branches'' of eigenvalues
$\lambda_i(h)$ that are convergent to zero as $h\to 0$ (the
``small'' eigenvalues) and the corresponding eigenspaces and
discuss the differentiable spectral sequence of the foliation,
which is a direct generalization of the Leray spectral sequence,
and its Hodge theoretic description.

We will consider two basic classes of foliations --- Riemannian
foliations and one-dimensional foliations defined by the orbits of
invariant flows on Riemannian Heisenberg manifolds.

We remark that the adiabatic limit is, up to scaling, an example of
collapsing (in general, without bounded curvature) in the sense of
\cite{Cheeger-Gromov}. For a discussion of the behavior of the
spectrum of the differential form Laplacian on a compact Riemannian
manifold under collapse, we refer, for instance, to
\cite{AmmannBar,Colbois,Fukaya87,jammes03,Lott02} and references
therein.

We are grateful to J. \'Alvarez L\'opez for useful discussions. We
also thank the referee for suggestions to improve the paper.

\section{Adiabatic limits and eigenvalue distribution}\label{s:adiab}

Let $(M,{\mathcal F})$ be a closed foliated manifold, endowed with
a Riemannian metric $g$. In this Section, we will discuss the
asymptotic behavior of the trace of $f(\Delta_h)$ in the adiabatic
limit.

\subsection{Riemannian foliations}
For Riemannian foliations, the problem was studied  in
\cite{adiab}. Recall (see, for instance,
\cite{Re1,Re,Molino82,Molino}) that a foliation ${\mathcal F}$ is
called Riemannian, if there exists a Riemannian metric $g$ on $M$
such that the induced metric $g_\tau$ on the normal bundle
$\tau=TM/F$ is holonomy invariant, or, equivalently, in any
foliated chart $\phi:U\to I^p\times I^q$ with local coordinates
$(x,y)$, the restriction $g_H$ of $g$ to $H=F^{\bot}$ is written
in the form $$ g_H=\sum_{\alpha,\beta=1}^q
g_{\alpha\beta}(y)\theta^{\alpha} \theta^{\beta}, $$ where
$\theta^{\alpha}\in H^*$ is the 1-form, corresponding to the form
$dy^\alpha$ under the isomorphism $H^*\cong T^*{\RR}^q$, and
$g_{\alpha\beta}(y)$ depend only on the transverse variables $y\in
\RR^q$. Such a Riemannian metric is called bundle-like.

It turns out that the adiabatic spectral limit on a Riemannian
foliation can be considered as a semiclassical spectral problem
for a Schr\"odinger operator on the leaf space $M/\cF$, and the
resulting asymptotic formula for the trace of $f(\Delta_h)$ can be
written in the form of the semiclassical Weyl formula for a
Schr\"odinger operator on a compact Riemannian manifold, if we
replace the classical objects, entering to this formula by their
noncommutative analogues. This observation provides a very natural
interpretation of the asymptotic formula for the trace of
$f(\Delta_h)$ (see Theorem~\ref{ad:main} below).

First, we transfer the operators $\Delta_{h}$ to the fixed Hilbert
space $L^2\Omega=L^{2}(M,\Lambda T^{*}M, g)$, using an isomorphism
$\Theta_h$ from $L^{2}(M,\Lambda T^{*}M , g_{h})$ to $L^2\Omega$
defined as follows. With respect to a bigrading on $\Lambda
T^{*}M$ given by
\[
\Lambda^k T^{*}M=\bigoplus_{i=0}^{k}\Lambda^{i,k-i}T^{*}M, \quad
\Lambda^{i,j}T^{*}M=\Lambda^{i}H^{*}\otimes \Lambda^{j}F^{*},
\]
we have, for $u \in L^{2}(M,\Lambda^{i,j}T^{*}M , g_{h})$,
\begin{equation}\label{e:theta}
\Theta_{h}u = h^{i}u.
\end{equation}
The operator $\Delta_h$ in $L^{2}(M,\Lambda T^{*}M , g_{h})$
corresponds under the isometry $\Theta_{h}$ to the operator
$L_{h}= \Theta_{h}\Delta_h\Theta_{h}^{-1}$ in $L^2\Omega$.

With respect to the above bigrading of $\Lambda T^*M$, the de Rham
differential $d$ can be written as
\[
d=d_F+d_H+\theta,
\]
where
\begin{enumerate}
\item  $ d_F=d_{0,1}: C^{\infty}(M,\Lambda^{i,j}T^*M)\to
C^{\infty}(M,\Lambda^{i,j+1}T^*M)$ is the tangential de Rham
differential, which is a first order tangentially elliptic
operator, independent of the choice of $g$;
\item $d_H=d_{1,0}: C^{\infty}(M,\Lambda^{i,j}T^*M)\to
C^{\infty}(M,\Lambda^{i+1,j}T^*M)$ is the transversal de Rham
differential, which is a first order transversally elliptic
operator;
\item $\theta=d_{2,-1}: C^{\infty}(M,\Lambda^{i,j}T^*M)\to
C^{\infty}(M,\Lambda^{i+2,j-1}T^*M)$ is a zeroth order
differential operator.
\end{enumerate}
One can show that
\[
d_h = \Theta_{h}d\Theta_{h}^{-1} = d_F + hd_H + h^{2}\theta,
\]
and the adjoint of $d_h$ in  $L^2\Omega$ is
\[
\delta_h=\Theta_{h}\delta \Theta_{h}^{-1} = \delta_F + h \delta_H
+ h^{2}\theta^{*}.
\]
Therefore, one has
\begin{equation*}
\begin{aligned}
L_h & = d_h\delta_h + \delta_h d_h \\ & =\Delta_F + h^2\Delta_H +
h^4\Delta_{\theta}+ hK_1+h^2K_2 +h^3K_3,
\end{aligned}
\end{equation*}
where $\Delta_F=d_Fd^*_F+d^*_Fd_F$ is the tangential Laplacian,
$\Delta_H=d_Hd^*_H+d^*_Hd_H$ is the transverse Laplacian,
$\Delta_{\theta}=\theta\theta^{*}+ \theta^{*}\theta$ and $K_2 =
d_F \theta^{*}+ \theta^{*} d_F + \delta_F\theta+\theta \delta_F$
are of zeroth order, and  $K_1 =  d_F \delta_H+ \delta_H d_F +
\delta_Fd_H+d_H \delta_F$ and $K_3 = d_H \theta^{*}+ \theta^{*}
d_H + \delta_H\theta+\theta \delta_H$ are first order differential
operators.

Suppose that ${\mathcal F}$ is a Riemannian foliation and $g$ is a
bundle-like metric. The key observation is that, in this case, the
transverse principal symbol of the operator $\delta_H$ is holonomy
invariant, and, therefore, the first order differential operator
$K_1$ is a {\em leafwise} differential operator. Using this fact,
one can show that the leading term in the asymptotic expansion of
the trace of $f(\Delta_h)$ or, equivalently, of the trace of
$f(L_h)$ as $h\to 0$ coincides with the leading term in the
asymptotic expansion of the trace of $f(\bar{L}_h)$ as $h\to 0$,
where
\[
\bar{L}_h=\Delta_F + h^2\Delta_H.
\]
More precisely, we have the following estimates (with some $C_1,
C_2>0$):
\[
|\tr f(L_h)|< C_1h^{-q}, \quad |\tr f(L_h)-\tr
f(\bar{L}_h)|<C_2h^{1-q}, \quad 0<h\leq 1,
\]
where we recall that $q$ denotes the codimension of $\cF$.

We observe that the operator $\bar{L}_h$ has the form of a
Schr\"odinger operator on the leaf space $M/\cF$, where $\Delta_H$
plays the role of the Laplace operator, and $\Delta_F$ the role of
the operator-valued potential on $M/\cF$.

Recall that, for a Schr\"odinger operator $H_h$ on a compact
Riemannian manifold $X$, $\dim X=n$, with a matrix-valued
potential $V\in C^\infty(X,{\mathcal L}(E))$, where $E$ is a
finite-dimensional Euclidean space and $V(x)^{*}=V(x)$:
\[
H_h=-h^2\Delta +V(x),\quad x\in X,
\]
the corresponding asymptotic formula (the semiclassical Weyl
formula) has the following form:
\begin{equation}\label{e:Weyl}
\operatorname{tr} f(H_h)=(2\pi)^{-n}h^{-n}\int_{T^*X}
\operatorname{Tr} f(p(x,\xi))\,dx\,d\xi+o(h^{-n}),\quad
h\rightarrow 0+,
\end{equation}
where $p\in C^\infty(T^*X,\cL(E))$
is the principal $h$-symbol of $H_h$:
\begin{equation*}
p(x,\xi)=|\xi|^2+V(x),\quad (x,\xi)\in T^{*}X.
\end{equation*}

Now we demonstrate how the asymptotic formula for the trace of
$f(\Delta_h)$ in the adiabatic limit can be written in a similar
form, using noncommutative geometry. (For the basic information on
noncommutative geometry of foliations, we refer the reader to
\cite{survey} and references therein.)

Let $G$ be the holonomy groupoid of $\cF$. Let us briefly recall
its definition. Denote by $\sim_h$ the equivalence relation on the
set of piecewise smooth leafwise paths $\gamma:[0,1] \rightarrow
M$, setting $\gamma_1\sim_h \gamma_2$ if $\gamma_1$ and $\gamma_2$
have the same initial and final points and the same holonomy maps.
The holonomy groupoid $G$ is the set of $\sim_h$ equivalence
classes of leafwise paths. $G$ is equipped with the source and the
range maps $s,r:G\rightarrow M$ defined by $s(\gamma)=\gamma(0)$
and $r(\gamma)=\gamma(1)$. Recall also that, for any $x\in M$, the
set $G^x=\{\gamma\in G : r(\gamma)=x\}$ is the covering of the
leaf $L_x$ through the point $x$, associated with the holonomy
group of the leaf. We will identify any $x\in M$ with the element
of $G$ given by the constant path $\gamma(t)=x, t\in [0,1]$.

Let $\lambda_L$ denote the Riemannian volume form on a leaf $L$
given by the induced metric, and $\lambda ^{x}, x\in M$, denote
the lift of $\lambda_{L_x}$ via the holonomy covering map
$s:G^x\to L_x$.

Denote by $\pi :N^*\cF\to M$ the conormal bundle to ${\cF}$ and by
${\cF}_N$ the linearized foliation in $N^*{\cF}$ (cf., for
instance, \cite{Molino,survey}). Recall that, for any $\gamma\in
G, s(\gamma)=x, r(\gamma)=y$, the codifferential of the
corresponding holonomy map defines a linear map $dh^*_{\gamma}:
N^*_y{\cF}\to N^*_x{\cF}$. Then the leaf of the foliation $\cF_N$
through $\nu\in N^*{\cF}$ is the set of all
$dh_{\gamma}^{*}(\nu)\in N^*{\cF}$, where $\gamma\in G,
r(\gamma)=\pi(\nu)$.

The holonomy groupoid $G_{{\mathcal F}_N}$ of the linearized
foliation ${\mathcal F}_N$ can be described as the set of all
$(\gamma,\nu)\in G\times {N}^*{\mathcal F}$ such that
$r(\gamma)=\pi(\nu)$. The source map $s_N:G_{{\mathcal
F}_N}\rightarrow {N}^*{\mathcal F}$ and the range map
$r_N:G_{{\mathcal F}_N}\rightarrow {N}^*{\mathcal F}$ are defined as
$s_N(\gamma,\nu)=dh_{\gamma}^{*}(\nu)$ and $r_N(\gamma,\nu)=\nu$. We
have a map $\pi_G:G_{{\mathcal F}_N}\rightarrow G$ given by
$\pi_G(\gamma,\nu)=\gamma$. Denote by ${\mathcal L}(\pi^*\Lambda
T^{*}M)$ the vector bundle on $G_{\cF_N}$, whose fiber at a point
$(\gamma,\nu)\in G_{\cF_N}$ is the space of linear maps
\[
(\pi^*\Lambda T^{*}M)_{s_N(\gamma,\nu)}\to (\pi^*\Lambda
T^{*}M)_{r_N(\gamma,\nu)}.
\]
There is a standard way (due to Connes \cite{Co79}) to introduce the
structure of involutive algebra on the space
$C^{\infty}_{c}(G_{{\mathcal F}_N}, {\mathcal L}(\pi^*\Lambda
T^{*}M))$ of smooth, compactly supported sections of ${\mathcal
L}(\pi^*\Lambda T^{*}M)$. For any $\nu \in N^*\cF$, this algebra has
a natural representation $R_\nu$ in the Hilbert space
$L^2(G^\nu_{\cF_N}, s^*_N(\pi^*\Lambda T^{*}M))$ that determines its
embedding to the $C^*$-algebra of all bounded operators in
$L^2(G_{\cF_N}, s^*_N(\pi^*\Lambda T^{*}M))$. Taking the closure of
the image of this embedding, we get a $C^*$-algebra $C^*(N^*\cF,
{\mathcal F}_N,\pi^{*}\Lambda T^{*}M)$, called the twisted foliation
$C^*$-algebra. The leaf space $N^*{\cF}/\cF_N$ can be informally
considered as the cotangent bundle to $M/\cF$, and the algebra
$C^*(N^*\cF, {\mathcal F}_N,\pi^{*}\Lambda T^{*}M)$ can be viewed as
a noncommutative analogue of the algebra of continuous vector-valued
differential forms on this singular space.

Let $g_N \in C^\infty(N^*\cF)$ be the fiberwise Riemannian metric
on $N^*{\mathcal F}$ induced by the metric on $M$. The principal
$h$-symbol of $\Delta_h$ is a tangentially elliptic operator in
$C^{\infty}(N^*{\mathcal F},\pi^{*} \Lambda T^{*}M)$ given by
\[
\sigma_h(\Delta_h) = \Delta_{{\mathcal F}_N}+g_N,
\]
where $\Delta_{{\mathcal F}_N}$ is the lift of the tangential
Laplacian $\Delta_F$ to a tangentially elliptic (relative to
${\mathcal F}_N$) operator in $C^{\infty}(N^*{\mathcal F},\pi^{*}
\Lambda T^{*}M)$, and $g_N$ denotes the multiplication operator in
$C^{\infty}(N^*{\mathcal F},\pi^{*} \Lambda T^{*}M)$ by the
function $g_N$. (Observe that $g_N$ coincides with the transversal
principal symbol of $\Delta_H$.) Consider $\sigma_h(\Delta_h)$ as
a family of elliptic operators along the leaves of the foliation
${\mathcal F}_N$ and lift these operators to the holonomy
coverings of the leaves. For any $\nu\in N^*\cF$, we get a
formally self-adjoint uniformly elliptic operator
$\sigma_h(\Delta_h)_\nu$ in $C^\infty(G^\nu_{\cF_N},
s^*_N(\pi^*\Lambda T^{*}M))$, which essentially self-adjoint in
the Hilbert space $L^2(G^\nu_{\cF_N}, s^*_N(\pi^*\Lambda
T^{*}M))$. For any $f\in S({\mathbb R})$, the family
$\{f(\sigma_h(\Delta_h)_\nu), \nu\in N^*\cF\}$ defines an element
$f(\sigma_h(\Delta_h))$ of the $C^*$-algebra $C^*(N^*\cF,
{\mathcal F}_N,\pi^{*}\Lambda T^{*}M)$.

The foliation $\cF_N$ has a natural transverse symplectic
structure, which can be described as follows. Consider a foliated
chart $\varkappa: U\subset M\rightarrow I^p\times I^q$ on $M$ with
coordinates $(x,y)\in I^p\times I^q$ ($I$ is the open interval
$(0,1)$) such that the restriction of $\cF$ to $U$ is given by the
sets $y={\rm const}$. One has the corresponding coordinate chart
in $T^*M$ with coordinates denoted by $(x,y,\xi,\eta)\in I^p\times
I^q\times \RR^p\times \RR^q$. In these coordinates, the
restriction of the conormal bundle $N^*{\mathcal  F}$ to $U$ is
given by the equation $\xi=0$. So we have a coordinate chart
$\varkappa_n : U_1\subset N^*{\mathcal  F} \longrightarrow
I^p\times I^q\times \RR^q$ on $N^*{\mathcal  F}$ with the
coordinates $(x,y,\eta)\in I^p\times I^q\times \RR^q$. Indeed, the
coordinate chart $\varkappa_n$ is a foliated coordinate chart for
${\mathcal  F}_N$, and the restriction of ${\mathcal F}_N$ to
$U_1$ is given by the level sets $y= {\rm const}, \eta={\rm
const}$. The transverse symplectic structure for $\cF_N$ is given
by the transverse two-form $\sum_jdy_j\wedge d\eta_j$.

The corresponding canonical transverse Liouville measure
$dy\,d\eta$ is holonomy invariant and, by noncommutative
integration theory \cite{Co79}, defines the trace
$\operatorname{tr}_{{\mathcal F}_N}$ on the $C^*$-algebra
$C^*(N^*\cF, {\mathcal F}_N,\pi^{*}\Lambda T^{*}M)$. Combining the
Riemannian volume forms $\lambda_L$ and the transverse Liouville
measure, we get a volume form $d\nu$ on $N^*\cF$. For any $k\in
C^{\infty}_{c}(G_{{\mathcal F}_N}, {\mathcal L}(\pi^*\Lambda
T^{*}M))$, its trace is given by the formula
\[
\operatorname{tr}_{{\mathcal F}_N}(k)=\int_{N^*\cF}k(\nu)d\nu.
\]
The trace $\operatorname{tr}_{{\mathcal F}_N}$ is a noncommutative
analogue of the integral over the leaf space $N^*{\cF}/\cF_N$ with
respect to the transverse Liouville measure. One can show that the
value of this trace on $f(\sigma_h(\Delta_h))$ is finite.

Replacing in the formula (\ref{e:Weyl}) the integration over
$T^{*}X$ and the matrix trace $\operatorname{Tr}$ by the trace
$\operatorname{tr}_{{\mathcal F}_N}$ and the principal $h$-symbol
$p$ by $\sigma_h(\Delta_h)$, we obtain the correct formula for
$\operatorname{tr} f(\Delta_{h})$ in the adiabatic limit.

\begin{theorem}[\cite{adiab}]
\label{ad:main} For any $f\in S({\mathbb R})$, the asymptotic
formula holds:
\begin{equation}\label{e:adiab}
\operatorname{tr} f(\Delta_{h}) =(2\pi)^{-q}h^{-q}
\operatorname{tr}_{{\mathcal F}_N} f(\sigma_h(\Delta_h))
+o(h^{-q}),\quad h\rightarrow 0.
\end{equation}
\end{theorem}

The formula (\ref{e:adiab}) can be rewritten in terms of the
spectral data of leafwise Laplace operators. We will formulate the
corresponding result for the spectrum distribution function
\[
N_h(\lambda)=\sharp \{i:\lambda_i(h)\leq \lambda\}.
\]

Restricting the tangential Laplace operator $\Delta_F$ to the
leaves of the foliation ${\mathcal F}$ and lifting the
restrictions to the holonomy coverings of leaves, we get the the
Laplacian $\Delta_x$ acting in $C^{\infty}_c(G^x,s^{*}\Lambda
T^{*}M)$. Using the assumption that ${\mathcal F}$ is Riemannian,
it can be checked that, for any $x\in M$, $\Delta_x$ is formally
self-adjoint in $L^2(G^x,s^{*}\Lambda T^{*}M)$, that, in turn,
implies its essential self-adjointness in this Hilbert space (with
initial domain $C^{\infty}_c(G^x,s^{*}\Lambda T^{*}M)$). For each
$\lambda \in {\mathbb R}$, let $E_x(\lambda)$ be the spectral
projection of $\Delta_x$, corresponding to the semi-axis $(-\infty
,\lambda ]$. The Schwartz kernels of the operators $E_x(\lambda)$
define a leafwise smooth section $e_\lambda$ of the bundle
$\cL(\Lambda T^{*}M)$ over $G$.

We introduce the spectrum distribution function $N_{\mathcal
F}(\lambda)$ of the operator $\Delta_F$ by the formula
\[
N_{\mathcal F}(\lambda)=\int_M \Tr e_\lambda(x)\,dx,\quad
\lambda\in {\mathbb R},
\]
where $dx$ denotes the Riemannian volume form on $M$. By
\cite{tang}, for any $\lambda  \in {\mathbb R}$, the function $\Tr
e_\lambda $ is a bounded measurable function on $M$, therefore,
the spectrum distribution function $N_{\mathcal F}(\lambda)$ is
well-defined and takes finite values.

As above, one can show that the family $\{E_x(\lambda): x\in M\}$
defines an element $E(\lambda)$ of the twisted von Neumann
foliation algebra $W^{*}(G,\Lambda T^{*}M)$, the holonomy
invariant transverse Riemannian volume form for $\cF$ defines a
trace $\tr_\cF$ on $W^{*}(G,\Lambda T^{*}M)$, and the right hand
side of the last formula can be interpreted as the value of this
trace on $E(\lambda)$.

\begin{theorem}[\cite{adiab}]
\label{intr} Let $(M,{\mathcal F})$ be a Riemannian foliation,
equipped with a bundle-like Riemannian metric $g$. Then the
asymptotic formula for $N_h(\lambda)$ has the following form:
\[
N_{h}(\lambda ) =h^{-q} \frac{(4\pi)^{-q/2}}{\Gamma((q/2)+1)}
\int_{-\infty}^{\lambda} (\lambda - \tau
)^{q/2}d_{\tau}N_{\mathcal F}(\tau )+o(h^{-q}),\quad h\rightarrow
0.
\]
\end{theorem}

\subsection{A linear foliation on the $2$-torus}\label{s:torus} In this
Section, we consider the simplest example of the situation studied
in the previous Section, namely, the example of a linear foliation
on the $2$-torus. So consider the two-di\-men\-si\-o\-nal torus
$\mathbb{T}^2=\mathbb{R}^2/\mathbb{Z}^2$ with the coordinates
$(x,y)\in \mathbb{R}^2$, taken modulo integer translations, and
the Euclidean metric $g$ on $\mathbb{T}^2$:
\[
g=d x^2+d y^2.
\]

Let $\widetilde{X}$ be the vector field on $\mathbb{R}^2$ given by
\[
\widetilde{X}=\frac{\partial}{\partial x}+\alpha
\frac{\partial}{\partial y},
\]
where $\alpha\in \mathbb R$. Since $\widetilde{X}$ is translation
invariant, it determines a vector field $X$ on $\mathbb{T}^2$. The
orbits of $X$ define a one-dimensional foliation $\mathcal F$ on
$\mathbb{T}^2$. The leaves of $\mathcal F$ are the images of the
parallel lines $\widetilde{L}_{(x_0,y_0)}=\{(x_{0}+t,
y_{0}+t\alpha) :t\in\mathbb{R}\}$, parameterized by $(x_0,y_0)\in
\mathbb R^2$, under the projection $\mathbb{R}^2 \rightarrow
\mathbb{T}^2$.

In the case when $\alpha$ is rational, all leaves of $\mathcal F$
are closed and are circles, and $\mathcal F$ is given by the
fibers of a fibration of $\mathbb{T}^2$ over $\mathbb{S}^1$. In
the case when $\alpha$ is irrational, all leaves of $\mathcal F$
are everywhere dense in $\mathbb{T}^2$.

The one-parameter family $g_{h}$ of Riemannian metrics on
$\mathbb{T}^2$ defined by (\ref{e:gh}) is given by
\[
g_h=
  \frac{1+h^{-2}\alpha^2}{1+\alpha^2} dx^2 + 2 \alpha
  \frac{1-h^{-2}}{1+\alpha^2}dx dy+
  \frac{\alpha^2+h^{-2}}{1+\alpha^2}dy^2.
\]
The Laplace operator (on functions) defined by $g_h$ has the form
$\Delta_h=\Delta_F+h^2\Delta_H$, where
 \[
\Delta_F=-\frac{1}{1+\alpha^2}\left( \frac{\partial}{\partial
x}+\alpha\frac{\partial}{\partial y} \right)^2, \quad
\Delta_H=-\frac{h^{2}}{1+\alpha^2}\left(
-\alpha\frac{\partial}{\partial x}+\frac{\partial}{\partial y}
\right)^2
\]
are the tangential and the transverse Laplace operators
respectively.

The operator $\Delta_h$ has a complete orthogonal system of
eigenfunctions
\[
u_{kl}(x,y)=e^{2\pi i(kx+ly)}, \quad (x,y)\in {\mathbb{T}^2},
\]
with the corresponding eigenvalues
\begin{equation}\label{e:eigen}
\lambda_{kl}(h)= (2\pi)^2 \left(\frac{1}{1+\alpha^2}(k+\alpha
l)^2+\frac{h^2}{1+\alpha^2}(-\alpha k+l)^2\right), \quad (k,l)\in
{\mathbb Z}^2.
\end{equation}

The eigenvalue distribution function of $\Delta_h$ has the form
\[
N_h= \# \{(k,l)\in {\mathbb Z}^2 : (2\pi)^2
\left(\frac{1}{1+\alpha^2}(k+\alpha
l)^2+\frac{h^2}{1+\alpha^2}(-\alpha k+l)^2\right) < \lambda\}.
 \]
Thus we come to the following problem of number theory:

\begin{prob}
Find the asymptotic for $h\rightarrow 0$ of the number of integer
points in the ellipse
\[
\{(\xi,\eta)\in \mathbb{R}^2: (2\pi)^2
\left(\frac{1}{1+\alpha^2}(\xi+\alpha
\eta)^2+\frac{h^2}{1+\alpha^2}(-\alpha
\xi+\eta)^2\right)<\lambda\}.
\]
\end{prob}

In the case when $\alpha$ is rational, this problem can be easily
solved by elementary methods of analysis. In the case when
$\alpha$ is irrational, such an elementary solution seems to be
unknown, and, in order to solve the problem, the connection of
this problem with the spectral theory of the Laplace operator and
with adiabatic limits plays an important role.

\begin{theorem}[\cite{torus-eng}]\label{th:main}
The following asymptotic formula for the spectrum distribution
function $N_h(\lambda)$ of the operator $\Delta_h$ for a fixed
$\lambda\in\mathbb{R}$ holds:
\medskip\par
1. For $\alpha\not\in\mathbb{Q},$
\begin{equation}\label{e:main1}
    N_{h}(\lambda ) =\frac{1}{4\pi}h^{-1}\lambda+o(h^{-1}), \quad
h\rightarrow 0.
\end{equation}

2. For $\alpha\in\mathbb{Q}$ of the form $\alpha=\frac{p}{q}$, where
$p$ and $q$ are coprime,
\begin{multline}\label{e:main2}
    N_{h}(\lambda ) \\ =h^{-1}
\sum_{\substack{k\in {\mathbb Z}\\
|k|<\frac{\sqrt{\lambda}}{2\pi}\sqrt{p^2+q^2}}}
\frac{1}{\pi\sqrt{p^2+q^2} }(\lambda -
\frac{4\pi^2}{p^2+q^2}k^2)^{1/2}+o(h^{-1}), \quad h\rightarrow 0.
\end{multline}
\end{theorem}

\begin{rem}
The asymptotic formulas (\ref{e:main1}) and (\ref{e:main2}) of
Theorem~\ref{th:main} look quite different. Nevertheless, it can
be shown that
\[
\lim_{\substack{p\to +\infty\\ q\to
+\infty}}\sum_{\substack{k\in {\mathbb Z}\\
|k|<\frac{\sqrt{\lambda}}{2\pi}\sqrt{p^2+q^2}}}
\frac{1}{\pi\sqrt{p^2+q^2} }(\lambda -
\frac{4\pi^2}{p^2+q^2}k^2)^{1/2}=\frac{1}{4\pi}\lambda.
\]
Indeed, one can write
\[
\sum_{k} \frac{1}{\pi\sqrt{p^2+q^2}}(\lambda -
\frac{4\pi^2}{p^2+q^2}k^2)^{1/2}=\sum_{k}g(\xi_k) \Delta\xi_k,
\]
where
\[
g(\xi)=\frac{1}{2\pi^2}(\lambda - \xi^2)^{1/2}, \quad
\xi_k=\frac{2\pi}{\sqrt{p^2+q^2}}k.
\]
This immediately implies that
\begin{multline*}
\lim_{\substack{p\to +\infty\\ q\to
+\infty}}\sum_{\substack{k\in {\mathbb Z}\\
|k|<\frac{\sqrt{\lambda}}{2\pi}\sqrt{p^2+q^2}}}
\frac{1}{\pi\sqrt{p^2+q^2} }(\lambda -
\frac{4\pi^2}{p^2+q^2}k^2)^{1/2} \\ =\frac{1}{2\pi^2}
\int_{-\sqrt{\lambda}}^{\sqrt{\lambda}}(\lambda - \xi^2)^{1/2}\,d\xi
=\frac{1}{4\pi}\lambda.
\end{multline*}
\end{rem}

We now show how to derive the asymptotic formulae of
Theorem~\ref{th:main} from Theorem~\ref{intr} (see
\cite{torus-eng} for more details).
\medskip\par
{\bf Case 1:} $\alpha\not\in\mathbb{Q}.$ In this case
$G=\mathbb{T}^2\times \mathbb{R}$. The source and the range maps
$s,r:G\rightarrow \mathbb{T}^2$ are defined for any
$\gamma=(x,y,t)\in G$ by $s(\gamma)=(x-t,y-\alpha t)$ and
$r(\gamma)=(x,y)$. For any $(x,y)\in \mathbb{T}^2$, the set
$G^{(x,y)}$ coincides with the leaf $L_{(x,y)}$ through $(x,y)$
and is diffeomorphic to $\mathbb R$:
\[
L_{(x,y)}=\{(x-t,y-\alpha t): t\in {\mathbb R}\}.
\]
The Riemannian volume form $\lambda^{(x,y)}$ on $L_{(x,y)}$ equals
$\sqrt{1+\alpha^2}\,dt$. Finally, the restriction of the operator
$\Delta_F$ to each leaf $L_{(x,y)}$ coincides with the operator
\[
A=-\frac{1}{1+\alpha^2}\frac{d^2}{dt^2},
\]
acting in the space $L^2({\mathbb{R}}, \sqrt{1+\alpha^2}\,dt)$.

Using the Fourier transform, it is easy to compute the Schwartz
kernel $E_\lambda(t,t_1)$ of the spectral projection
$\chi_\lambda(A)$ of the operator $A$, corresponding to the
semi-axis $(-\infty,\lambda]$ (relative to the volume form
$\sqrt{1+\alpha^2}\,dt$):
\[
E_\lambda(t,t_1)=\frac{1}{2\pi\sqrt{1+\alpha^2}} \int_{\mathbb{R}}
e^{i(t-t_1)\xi}
\chi_\lambda\left(\frac{|\xi|^2}{1+\alpha^2}\right)d\xi.
\]
Then, for any $\gamma=(x,y,t)\in G=\mathbb{T}^2\times \mathbb{R}$,
we have
\[
e_\lambda(\gamma)=E_{\lambda}(0,t)=\frac{1}{2\pi\sqrt{1+\alpha^2}}
\int_{\mathbb{R}}{e^{-it\xi}\chi_\lambda\left(\frac{|\xi|^2}{1+\alpha^2}\right)d\xi}.
\]
The restriction of $e_\lambda$ to $\mathbb{T}^2$ is given by
\[
e_\lambda(x,y)=E_{\lambda}(0,0)=\frac{1}{2\pi\sqrt{1+\alpha^2}}
\int_{\mathbb{R}}\chi_\lambda\left(\frac{|\xi|^2}{1+\alpha^2}\right)d\xi
= \frac{1}{\pi}\sqrt{\lambda}, \quad \lambda>0.
\]
We get that the spectrum distribution function $N_{\mathcal
F}(\lambda )$ of the operator $\Delta_{F}$ has the form:
\[
N_\mathcal{F}(\lambda)=\int_{\mathbb{T}^2}{e_\lambda(x,y)dxdy}=
 \frac{1}{\pi}\sqrt{\lambda}, \quad \lambda>0.
 \]
By Theorem~\ref{intr}, we obtain
\begin{align*}
N_{h}(\lambda ) & =h^{-1} \frac{1}{\pi} \int_{-\infty}^{\lambda}
(\lambda - \tau )^{1/2}d_{\tau}N_{\mathcal F}(\tau )+o(h^{-1})\\
&=\frac{1}{4\pi}h^{-1}\lambda+o(h^{-1}), \quad h\rightarrow 0.
\end{align*}
\medskip\par
{\bf Case 2:} $\alpha\in\mathbb{Q}$ of the form
$\alpha=\frac{p}{q}$, where $p$ and $q$ are coprime. In this case,
the holonomy groupoid is $\mathbb{T}^2\times
(\mathbb{R}/{q\mathbb{Z}})$. The leaf $L_{(x,y)}$ through any
$(x,y)$ is the circle $\{(x+t,y+\alpha t): t\in
\mathbb{R}/{q\mathbb{Z}}\}$ of length $l=\sqrt{p^2+q^2}$. The
restriction of the operator $\Delta_F$ to each $L_{(x,y)}$
coincides with the operator
\[
A=-\frac{1}{1+\alpha^2}\frac{d^2}{dt^2},
\]
acting in the space $L^2(\mathbb{R}/q\mathbb{Z},
\sqrt{1+\alpha^2}\,dt)$.

Using the Fourier transform, it is easy to see that the kernel of
the spectral projection $\chi_\lambda(A)$ in
$L^2(\mathbb{R}/q\mathbb{Z}, \sqrt{1+\alpha^2}\,dt)$ is given by
the formula
\[
E_{\lambda}(t,t_1)=\frac{1}{\sqrt{p^2+q^2}}
\sum_{\substack{k\in {\mathbb Z}\\
|k|<\frac{\sqrt{\lambda}}{2\pi}\sqrt{p^2+q^2}}} {e^{\frac{2\pi
i}{q}k(t-t_1)}}.
    \]
For any $\gamma=(x,y,t)\in G=\mathbb{T}^2
\times(\mathbb{R}/q\mathbb{Z})$, we have
\[
e_\lambda(\gamma)=E_\lambda(0,t)=\frac{1}{\sqrt{p^2+q^2}}
\sum_{\substack{k\in {\mathbb Z}\\
|k|<\frac{\sqrt{\lambda}}{2\pi}\sqrt{p^2+q^2}}}{e^{-\frac{2\pi
i}{q}k t}}.
\]
We get that the spectrum distribution function $N_{\mathcal
F}(\lambda )$ of $\Delta_{F}$ is of the form:
\[
N_{\mathcal{F}}(\lambda )= \int_{\mathbb{T}^2}{e_\lambda(x,y)dxdy}
=\frac{1}{\sqrt{p^2+q^2}}\#\{k\in {\mathbb Z} :
|k|<\frac{\sqrt{\lambda}}{2\pi}\sqrt{p^2+q^2}\}.
\]
By Theorem~\ref{intr}, we obtain for $h\rightarrow 0$
\[
\begin{split}
 N_{h}(\lambda )
& =h^{-1} \frac{1}{\pi} \int_{-\infty}^{\lambda} (\lambda - \tau
)^{1/2}d_{\tau}N_{\mathcal F}(\tau )+o(h^{-1})\\
& =h^{-1} \frac{1}{\pi\sqrt{p^2+q^2} }
\sum_{\substack{k\in {\mathbb Z}\\
|k|<\frac{\sqrt{\lambda}}{2\pi}\sqrt{p^2+q^2}}} (\lambda -
\frac{4\pi^2}{p^2+q^2}k^2)^{1/2}+o(h^{-1}).
\end{split}
\]

\subsection{Riemannian Heisenberg manifolds}\label{s:Heis0}
In this Section we consider the adiabatic limits associated with
one-dimensional foliations given by the orbits of invariant flows
on Riemannian Heisenberg manifolds. These foliations are examples
of non-Riemannian foliations.

Recall that the real three-dimensional Heisenberg group $H$ is the
Lie subgroup of $\operatorname{GL}(3,\mathbb{R})$ consisting of
all matrices of the form
\[
\gamma (x,y,z)=\begin{bmatrix}
 {1} & {x} & {z} \\
 {0} & {1} & {y}  \\
 0 & 0 & {1} \\
\end{bmatrix}, \quad x,y,z\in\mathbb{R}.
\]
Its Lie algebra $\mathfrak{h}$ is the Lie subalgebra of
$gl(3,\mathbb{R})$ consisting of all matrices of the form
\[
X(x,y,z)=\begin{bmatrix}
{0} & {x} & {z} \\
{0} & {0} & {y}  \\
 0 & 0 & {0} \\
\end{bmatrix}, \quad x,y,z\in\mathbb{R}.
\]

A Riemannian Heisenberg manifold $M$ is defined to be a pair
$(\Gamma \backslash H,g)$, where $\Gamma =\{ \gamma (x,y,z): x,y,z
\in \mathbb Z \}$ is a uniform discrete subgroup of $H$ and $g$ is
a Riemannian metric on $\Gamma \backslash H$ whose lift to $H$ is
left $H$-invariant.

It is easy to see that $g$ is uniquely determined by the value of
its lift to $H$ at the identity $\gamma(0,0,0)$, that is, by a
symmetric positive definite $3\times 3$-matrix.

In the following, we will assume that the metric $g$ corresponds to
a $3\times 3$-matrix of the form
\begin{equation}\label{e:g1}
\begin{pmatrix}
  h_{11} & h_{12} & 0 \\
  h_{12} & h_{22} & 0 \\
  0 & 0 & g_{33} \\
\end{pmatrix}.
\end{equation}

The lift of $g$ to $H$ is given by the formula
\begin{multline*}
g(\gamma(x,y,z))=
h_{11}dx^2+2h_{12}dx\,dy+h_{22}dy^2+g_{33}(dz-x\,dy)^2,
\\ (x,y,z) \in {\mathbb R}^3.
\end{multline*}
The corresponding Laplace operator has the form
  \begin{multline*}
\Delta =-\left\{ \frac{1}{h_{11}h_{22}-h^2_{12}} \right. \left[
h_{22}\frac{\partial^{2}}{\partial x^{2}} -
h_{12}\right.\left[\frac{\partial}{\partial x}\right.
\left(\frac{\partial }{\partial y} + x\frac{\partial }{\partial z}
\right)\\  + \left(\frac{\partial }{\partial y} + x\frac{\partial
}{\partial z}\right)\left. \frac{\partial}{\partial x}\right]
\left. +h_{11}\left(\frac{\partial }{\partial y} + x\frac{\partial
}{\partial z}\right)^{2} \right] + \left.
\frac{1}{g_{33}}\frac{\partial^{2}}{\partial z^{2}} \right\}.
\end{multline*}

\begin{thm}[\cite{Gordon-Wilson}]\label{thm:0spectrum}
The spectrum of the Laplace operator $\Delta $ on functions on $M$
(with multiplicities) has the form
\[
{\rm spec}\,\Delta =\Sigma_{1}\cup \Sigma_{2},
\]
where
\begin{multline*}
\begin{aligned}
\Sigma_{1} & = \{\lambda(a,b)
=4\pi^{2}\frac{h_{22}a^2-2h_{12}ab+h_{11}b^2}{h_{11}h_{22}-h^2_{12}}
:a,b\in \mathbb Z\},\\
\Sigma_{2} & = \{\mu(c,k)=\frac{4\pi^2 c^{2}}{g_{33}}+\frac{2\pi c
(2k+1)}{\sqrt{h_{11}h_{22}-h^2_{12}}}\ \text{with mult.}\ 2c :
\end{aligned}
\\
c\in \mathbb Z^{+},\quad k\in \mathbb Z^{+}\cup\{0\}\}.
\end{multline*}
\end{thm}

\begin{rem}
As shown in \cite{Gordon-Wilson}, for an arbitrary left
$H$-invariant metric $g$ on $H$, there exists a left $H$-invariant
metric $g_1$, which corresponds to a $3\times 3$-matrix of the
form (\ref{e:g1}), such that Riemannian Heisenberg manifolds
$(\Gamma \backslash H,g)$ and $(\Gamma \backslash H,g_1)$ are
isometric. Therefore, Theorem~\ref{thm:0spectrum} provides a
solution of the problem of calculation of the spectrum of the
Laplace operator on functions for an arbitrary Riemannian
Heisenberg manifold.
\end{rem}

Now we assume that the metric $g$ on $M$ corresponds to a $3\times
3$-matrix of the form
\[
\begin{pmatrix}
  h_{11} & 0 & 0 \\
  0 & h_{22} & 0 \\
  0 & 0 & g_{33} \\
\end{pmatrix}.
\]
In this case, one can write down explicitly all the eigenfunctions
of the corresponding Laplace operator on functions. This fact plays
an important role in the proof of the following theorem.

\begin{thm}[\cite{Heisenberg}]\label{thm:1spectrum}
The spectrum of the Laplace operator $\Delta $ on differential one
forms on $M$ (with multiplicities) has the form
\[
{\rm spec}\,\Delta =\Sigma_{1}\cup \Sigma_{2}\cup \Sigma_{3},
\]
where
\begin{multline*}
\Sigma_{1}  = \{\lambda_\pm(a,b)
=4\pi^{2}\left(\frac{a^2}{h_{11}}+\frac{b^2}{h_{22}}\right) \\ +
\frac{\frac{g_{33}}{h_{11}h_{22}}\pm
\sqrt{\frac{g^2_{33}}{h^2_{11}h^2_{22}}+
16\pi^{2}\frac{g_{33}}{h_{11}h_{22}}(\frac{a^2}{h_{11}}+\frac{b^2}{h_{22}})}}{2}
\ \text{with mult.}\ 2 :a,b\in \mathbb Z\},
\end{multline*}

\begin{multline*}
\Sigma_{2} = \{\mu(c,k)=\frac{4\pi^2 c^{2}}{g_{33}}+\frac{2\pi c
(2k+1)}{\sqrt{h_{11}h_{22}}}\ \text{with mult.}\ 2c :
\\
c\in \mathbb Z^{+}, k\in \mathbb Z^{+}\cup\{0\}\}.
\end{multline*}

\begin{multline*}
\Sigma_{3}  = \{\mu_\pm(c,k)= \frac{4\pi^2c^2}{g_{33}}+\frac{2\pi
c (2k+1)}{\sqrt{h_{11}h_{22}}}
\\ + \frac{\frac{g_{33}}{h_{11}h_{22}} \pm
\sqrt{(\frac{4\pi c
}{\sqrt{h_{11}h_{22}}}+\frac{g_{33}}{h_{11}h_{22}})^2
+8k\frac{2\pi c g_{33}}{(\sqrt{h_{11}h_{22}})^3}}}{2}\ \text{with
mult.}\ 2c:
\\
c\in \mathbb Z^{+}, k\in \mathbb Z^{+}\cup\{0\}\}.
\end{multline*}
\end{thm}

We refer the reader to \cite{AmmannBar} for a similar calculation
of the spectrum of the Dirac operator on Riemannian Heisenberg
manifolds.

Let $\alpha  \in \mathbb{R}$. Consider the left-invariant vector
field on $H$ associated with
\[
  X(1,\alpha,0)=\begin{bmatrix}
{0} & {1} & {0} \\
{0} & {0} & {\alpha}  \\
 0 & 0 & {0} \\
\end{bmatrix} \in \mathfrak{h}.
\]
Since $X(1,\alpha,0)$ is a left-invariant vector field, it
determines a vector field on $M=\Gamma \backslash H$. The orbits
of this vector field define a one-dimensional foliation $\mathcal
F$ on $M$. The leaf through a point $\Gamma \gamma(x,y,z)\in M$ is
described as
\[
L_{\Gamma \gamma(x,y,z)} =\{ \Gamma \gamma(x+t,y+\alpha t,z+\alpha
t x+\frac{\alpha t^2}{2})\in \Gamma \backslash H: t\in {\mathbb
R}\}.
\]

We assume that $g$ corresponds to the identity $3\times 3$-matrix.
Consider the adiabatic limit associated with the Riemannian
Heisenberg ma\-ni\-fold $(\Gamma \backslash H,g)$ and the
one-dimensional foliation $\mathcal F$. The Riemannian metric
$g_h$ on $\Gamma \backslash H$ defined by (\ref{e:gh}) corresponds
to the matrix
\[
 \begin{pmatrix}
  \frac{1+h^{-2}\alpha^2}{1+\alpha^2}&\alpha\frac{1-h^{-2}}{1+\alpha^2}&0\\
  \alpha\frac{1-h^{-2}}{1+\alpha^2}&\frac{\alpha^2+h^{-2}}{1+\alpha^2}&0\\
  0&0&h^{-2}
 \end{pmatrix}, \quad h>0.
\]
The corresponding Laplacian (on functions) on the group $H$ has
the form:
\begin{multline*}
\Delta_h=-\frac{1}{1+\alpha^2}\Bigg[   \left(
\frac{\partial}{\partial x}+\alpha\left(\frac{\partial}{\partial
y}+x\frac{\partial}{\partial z} \right) \right)^2 \\ +h^2\left(
-\alpha\frac{\partial}{\partial x}+\frac{\partial}{\partial
y}+x\frac{\partial}{\partial z} \right)^2
\Bigg]-h^2\frac{\partial^2}{\partial z^2}.
\end{multline*}

Using an explicit computation of the heat kernel on the Heisenberg
group, one can show the following asymptotic formula.

\begin{thm}[\cite{matsb}] \label{t:trace}
For any $t>0$, we have as $h\to 0$
\begin{equation}\label{e:nil}
\operatorname{tr} e^{-t\Delta_h}=
\frac{h^{-2}}{4\pi}\int_{-\infty}^{+\infty}
\frac{\eta}{\sinh(t\eta)} e^{-t\eta^2} d\eta + o(h^{-2}).
\end{equation}
\end{thm}

\begin{rem}
The formula (\ref{e:nil}) looks quite different from what we have in
the case of a Riemannian foliation. For instance, if $\cF$ is a
one-dimensional Riemannian foliation on a three-dimensional closed
Riemannian manifold $M$ given by the orbits of an non-singular
isometric flow such that the set of closed orbits has measure zero,
then, by Theorem~\ref{ad:main} (or, equivalently, by
Theorem~\ref{intr}), the asymptotic formula for the trace of the
heat operator $e^{-t\Delta_h}$ in the adiabatic limit has the
following form: for any $t>0$,
\begin{align*}
\operatorname{tr} e^{-t\Delta_h} & = \frac{h^{-2}}{4\pi t
}\int_{-\infty}^{+\infty} e^{-t\eta^2} d\eta + o(h^{-2})\\ & =
\frac{h^{-2}}{4\sqrt{\pi t^3}} + o(h^{-2}), \quad h\to 0.
\end{align*}
So, in comparison with the case of Riemannian foliations, the
formula (\ref{e:nil}) contains an additional factor
$\frac{t\eta}{\sinh(t\eta)}$ related with the distortion of the
transverse part of the Riemannian metric along the orbits of the
flow.
\end{rem}

\begin{rem}
It would be quite interesting to write the formula (\ref{e:nil}) in
a form similar to the formula (\ref{e:adiab}).
\end{rem}

\section{Adiabatic limits and differentiable spectral
sequence}\label{s:spec} In this Section, we will discuss the
problem of ``small eigenvalues'' of the Laplace operator in the
adiabatic limit and its relation with the differentiable spectral
sequence of the foliation. We will start with some background
information on the differentiable spectral sequence.

\subsection{Preliminaries on the differentiable spectral sequence}
As usual, let $\cF$ be a codimension $q$ foliation on a closed
manifold $M$. The differentiable spectral sequence $(E_k,d_k)$ of
$\cF$ is a direct generalization of (the differentiable version
of) the Leray spectral sequence for fibrations, which converges to
the de~Rham cohomology of $M$.

Denote by $\Omega$ the space of smooth differential forms and by
$\Omega^r$ the space of smooth differential $r$-forms on $M$.
Similar to the bundle case, the {\em differentiable spectral
sequence\/} $(E_k,d_k)$ of $\cF$ is defined by the decreasing
filtration by differential subspaces
$$\Omega=\Omega_0\supset\Omega_1\supset\cdots
\supset\Omega_q\supset\Omega_{q+1}=0\;,$$ where the space
$\Omega_k^r$ of $r$-forms of filtration degree $\geq k$ consists
of all $\omega\in\Omega^r$ such that $i_X\omega=0$ for all
$X=X_1\wedge\cdots\wedge X_{r-k+1}$, where the $X_i$ are vector
fields tangent to the leaves. Roughly speaking, $\omega$ in
$\Omega_k^r$ iff it is of degree $\ge k$ transversely to the
leaves.

Recall that the induced spectral sequence $(E_k,d_k)$ is defined
in the following standard way (see, for instance,
\cite{McCleary}):
\begin{alignat*}{3}
Z_k^{u,v}&=\Omega_u^{u+v}\cap
d^{-1}\left(\Omega_{u+k}^{u+v+1}\right)\;, &\quad
Z_\infty^{u,v}&=\Omega_u^{u+v}\cap \ker d\;,\\
B_k^{u,v}&=\Omega_u^{u+v}\cap
d\left(\Omega_{u-k}^{u+v-1}\right)\;, &\quad
B_\infty^{u,v}&=\Omega_u^{u+v}\cap \im d\;,\\
E_k^{u,v}&=\frac{Z_k^{u,v}}{Z_{k-1}^{u+1,v-1}+B_{k-1}^{u,v}}\;,
&\quad
E_\infty^{u,v}&=\frac{Z_\infty^{u,v}}{Z_\infty^{u+1,v-1}+B_\infty^{u,v}}\;.
\end{alignat*}
We assume $B_{-1}^{u,v}=0$, so
$E_0^{u,v}=\Omega_u^{u+v}/\Omega_{u+1}^{u+v}$. Each homomorphism
$d_k:E_k^{u,v}\to E_k^{u+k,v-k+1}$ is canonically induced by $d$.

The terms $E_1^{0,\ast}$ and $E_2^{\ast,0}$ are respectively
called {\em leafwise cohomology\/} and {\em basic cohomology\/},
and $E_2^{\ast,p}$ is isomorphic to the {\em transverse
cohomology\/} \cite{Haefliger80} (also called {\em Haefliger
cohomology\/}).

The $C^\infty$~topology of $\Omega$ induces a topological vector
space structure on each term $E_k$ such that $d_k$ is continuous.
A subtle problem here is that $E_k$ may not be Hausdorff
\cite{Haefliger80}. So it makes sense to consider the subcomplex
given by the closure of the trivial subspace, $\bar0_k\subset
E_k$, as well as the quotient complex ${\widehat
E}_k=E_k/\bar0_k$, whose differential operator will be also
denoted by $d_k$.

\subsection{Riemannian foliations}
For a Riemannian foliation $\cF$, each term $E_k$ of the
differentiable spectral sequence $(E_k,d_k)$ is Hausdorff of finite
dimension if $k\geq2$, and $H(\bar0_1)=0$. So $E_k\cong{\widehat
E}_k$ for $k\geq2$. The proof of this result given in \cite{Masa92}
uses the structure theorem for Riemannian foliations due to Molino
\cite{Molino82,Molino} to reduce the problem to transitive
foliations, and, for transitive foliations, it uses a construction
of a parametrix for the de Rham complex given by Sarkaria
\cite{Sarkaria}. Moreover, it turns out that, for $k\geq2$, the
terms $E_k$ are homotopy invariants of Riemannian foliations
\cite{AlvMasa1}. (This result generalizes a previous work showing
the topological invariance of the basic cohomology
\cite{KacimiNicolau93}.)

Now return to adiabatic limits. So let $g$ be a Riemannian metric
on $M$ and $g_h$ be the one-parameter family of metrics defined by
(\ref{e:gh}). Denote by $\Delta^r_h$ the Laplace operator on
differential $r$-forms on $M$ defined by $g_h$, and by
\[
0\leq\lambda_0^r(h)\leq\lambda_1^r(h)\leq\lambda_2^r(h)\leq\cdots
\]
its spectrum (with multiplicities). It is well known that the
eigenvalues of the Laplacian on differential forms vary
continuously under continuous perturbations of the metric, and
thus the ``branches'' of eigenvalues $\lambda_i^r(h)$ depend
continuously on $h>0$. In this Section, we shall only consider the
``branches'' $\lambda_i^r(h)$ that are convergent to zero as $h\to
0$; roughly speaking, the ``small'' eigenvalues. The asymptotics
as $h\to 0$ of these metric invariants are related to the
differential invariant $\widehat E_1^r$ and the homotopy
invariants $E_k^r$, $k\geq2$, as follows.

\begin{theorem}[\cite{seq}]\label{small eigenvalues} With the above notation,
for Riemannian foliations on closed Riemannian manifolds we have
\begin{align*}
\dim{\widehat E}_1^r&= \sharp\,\left\{i\ \left|\ \lambda_i^r(h)=
O\left(h^2\right) \quad\text{as}\quad h\to 0\right.\right\}\;,\\
\dim E_k^r&= \sharp\,\left\{i\ \left|\ \lambda_i^r(h)=
O\left(h^{2k}\right) \quad\text{as}\quad h\to
0\right.\right\}\;,\quad k\geq2\;.
\end{align*}
\end{theorem}

We refer to \cite{jammes05} for a particular form of this theorem
in the case of Riemannian flows.

As a part of the proof of Theorem~\ref{small eigenvalues} and also
because of its own interest, the asymptotics of eigenforms of
$\Delta_{h}$ corresponding to ``small'' eigenvalues were also
studied. This study was begun in \cite{MazzeoMelrose} for the case
of Riemannian bundles, and continued in \cite{Forman95} for
general complementary distributions.

Here we formulate the results obtained in \cite{seq} for the case
of Riemannian foliations. Recall that $\Theta_h$ is an isomorphism
of Hilbert spaces, which moves our setting to the fixed Hilbert
space $L^2\Omega$ (see (\ref{e:theta})). The ``rescaled
Laplacian'' $L_h=\Theta_h\Delta_{h}\Theta_h^{-1}$ has the same
spectrum as $\Delta_{h}$, and eigenspaces of $\Delta_{h}$ are
transformed into eigenspaces of $L_h$ by $\Theta_h$. It turns out
that eigenspaces of $L_h$ corresponding to ``small'' eigenvalues
are convergent as $h\to0$ when the metric $g$ is bundle-like, and
the limit is given by a nested sequence of bigraded subspaces,
$$
\Omega\supset\cH_1\supset\cH_2\supset\cH_3\supset\cdots\supset\cH_\infty\;.
$$ The definition of $\cH_1,\cH_2$ was given in \cite{AlvKordy1}
as a Hodge theoretic approach to $(E_1,d_1)$ and $(E_2,d_2)$,
which is based on the study of leafwise heat flow. The space
${\cH}_1$ is defined as the space of smooth leafwise harmonic
forms:
\[
{\cH}_1=\{\omega\in \Omega : \Delta_F\omega=0 \}.
\]
As shown in \cite{AlvKordy1}, the orthogonal projection in
$L^2\Omega$ on the kernel of $\Delta_F$ in $L^2\Omega$ restricts
to smooth differential forms, yielding an operator $\Pi : \Omega
\to {\cH}_1$. We define the operator $d_1$ on ${\cH}_1$ as
$d_1=\Pi d_H$. The adjoint of $d_1$ in $\cH_1$ equals
$\delta_1=\Pi\delta_H$. Finally, we take
$\Delta_1=d_1\delta_1+\delta_1d_1$ on ${\cH}_1$ and put
\[
\cH_2=\ker\Delta_1.
\]

The other spaces $\cH_k$ are defined in \cite{seq} as an extension
of this Hodge theoretic approach to the whole spectral sequence
$(E_k,d_k)$. In particular,
\[
\cH_1\cong\widehat E_1\;,\quad\cH_k\cong E_k\;,\quad
k=2,3,\ldots,\infty\;,
\]
as bigraded topological vector spaces. Thus this sequence
stabilizes (that is, $\cH_k=\cH_\infty$ for $k$ large enough)
because the differentiable spectral sequence is convergent in a
finite number of steps. The convergence of eigenforms
corresponding to ``small'' eigenvalues is precisely stated in the
following result, where $L^2\cH_1$ denotes the closure of $\cH_1$
in $L^2\Omega$.

\begin{theorem} For any Riemannian foliation on a closed manifold
with a bundle-like metric, let $\omega_i$ be a sequence in
$\Omega^r$ such that $\|\omega_i\|=1$ and
\[
\left\langle L_{h_i}\omega_i,\omega_i\right\rangle\in
o\left(h_i^{2(k-1)}\right)
\]
for some fixed integer $k\ge1$ and some sequence $h_i\to0$. Then
some subsequence of the $\omega_i$ is strongly convergent, and its
limit is in $L^2\cH_1^r$ if $k=1$, and in $\cH_k^r$ if $k\geq2$.
\end{theorem}

To simplify notation let $m_1^r=\dim\widehat E_1^r$, and let
$m_k^r=\dim E_k^r$ for each $k=2,3,\ldots,\infty$. Thus
Theorem~\ref{small eigenvalues} establishes $\lambda_i^r(h)=
O\left(h^{2k}\right)$ for $i\leq m_k^r$, yielding
$\lambda_i^r(h)\equiv0$ for $i$ large enough. For every $h>0$,
consider the nested sequence of graded subspaces
$$\Omega\supset\cH_1(h)\supset\cH_2(h)\supset\cH_3(h)\supset\cdots\supset\cH_\infty(h)\;,$$
where $\cH_k^r(h)$ is the space generated by the eigenforms of
$\Delta_h$ corresponding to eigenvalues $\lambda_i^r(h)$ with
$i\leq m_k^r$; in particular, we have
$\cH_k(h)=\cH_\infty(h)=\ker\Delta_h$ for $k$ large enough. Set
also $\cH_k(0)=\cH_k$. We have $\dim\cH_k^r(h)=m_k^r$ for all
$h>0$, so the following result is a sharpening of
Theorem~\ref{small eigenvalues}.

\begin{cor}\label{c:seq} For any Riemannian foliation on a closed manifold with
a bundle-like metric and $k=2,3,\ldots,\infty$, the assignment
$h\mapsto\cH_k^r(h)$ defines a continuous map from $[0,\infty)$ to
the space of finite dimensional linear subspaces of $L^2\Omega^r$
for all $r\geq0$. If $\dim\widehat E_1^r<\infty$, then this also
holds for $k=1$. \end{cor}

By the standard perturbation theory, the map $h\mapsto\cH_k^r(h)$
is, clearly, $C^\infty$ on $(0,\infty)$ for any Riemannian foliation
on a closed manifold with a bundle-like metric,
$k=2,3,\ldots,\infty$ and $r\geq0$. As shown in
\cite{MazzeoMelrose}, this map is $C^\infty$ up to $h=0$, if the
foliation is given by the fibers of a Riemannian fibration. In the
next section, we will see an example of a Riemannian foliation and a
bundle-like metric such that the map $h\mapsto\cH_k^r(h)$ is not
$C^\infty$ at $h=0$.

\subsection{A linear foliation on the $2$-torus}
In this Section, we consider the simplest example of the situation
studied in the previous section, namely --- the example of a linear
foliation on the $2$-torus. So, as in Section~\ref{s:torus},
consider the two-dimensional torus
$\mathbb{T}^2=\mathbb{R}^2/\mathbb{Z}^2$ with the coordinates
$(x,y)$, the one-dimensional foliation $\cF$ defined by the orbits
of the vector field $\widetilde{X}=\frac{\partial}{\partial
x}+\alpha \frac{\partial}{\partial y}$, where $\alpha\in \mathbb R$,
and the Euclidean metric $g=d x^2+d y^2$ on $\mathbb{T}^2$. The
eigenvalues of the corresponding Laplace operator $\Delta_h$
(counted with multiplicities) are described as follows:
\begin{gather*}
{\rm spec}\,\Delta^0_h={\rm spec}\,\Delta^2_h= \{ \lambda_{kl}(h) :
(k,l)\in {\mathbb Z}^2\}, \\ {\rm spec}\,\Delta^1_h = \{
\lambda_{k_1l_1}(h)+\lambda_{k_2l_2}(h) : (k_1,l_1)\in {\mathbb
Z}^2, (k_2,l_2)\in {\mathbb Z}^2\},
\end{gather*}
where $\lambda_{kl}(h)$ are given by (\ref{e:eigen}). So, for
$\alpha\notin \QQ$, small eigenvalues appear only if $(k,l)=(0,0)$
and $(k_1,l_1)=(k_2,l_2)=(0,0)$ and have the form
\begin{equation}\label{e:form}
\lambda^0_0(h)=\lambda^2_0(h)=0, \quad
\lambda^1_0(h)=\lambda^1_1(h)=0.
\end{equation}
For $\alpha\in\mathbb{Q}$ of the form $\alpha=\frac{p}{q}$, where
$p$ and $q$ are coprime, small eigenvalues appear only if
$(k,l)=t(p,q), t\in \ZZ,$ and $(k_1,l_1)=t_1(p,q),
(k_2,l_2)=t_2(p,q), t_1, t_2\in \ZZ$. So there are infinitely many
different branches of eigenvalues $\lambda_h$ with
$\lambda_h=O(h^2)$ as $h\to0$, and all the branches of eigenvalues
$\lambda_h$ with $\lambda_h=O(h^4)$ as $h\to0$ are given by
(\ref{e:form}).

Now let us turn to the differential spectral sequence. By a
straightforward computation, one can show that
\[
E^{u,v}_2=E^{u,v}_\infty={\mathbb R}, \quad u=0,1, \quad v=0,1,
\]
that agrees with the above description of small eigenvalues.

The case of $E_1$ is more interesting. First of all, it depends on
whether $\alpha$ is rational or not. For $\alpha\in \QQ$, $\cF$ is
given by the fibers of a trivial fibration $\mathbb{T}^2\to S^1$,
and, therefore, for any $u=0,1$ and $v=0,1$, we have
\[
E^{u,v}_1={\widehat E}^{u,v}_1 = \Omega^u(S^1)\otimes
H^v(S^1)=C^\infty(S^1).
\]
For $\alpha\notin \QQ$, we have
\[
{\widehat E}^{u,v}_1={\mathbb R}, \quad u=0,1, \quad v=0,1.
\]
The description of $E_1$ is more complicated and depends on the
diophantine properties of $\alpha$. Recall that $\alpha\notin \QQ$
is called diophantine, if there exist $c>0$ and $d>1$ such that,
for any $p\in \ZZ\setminus\{0\}$ and $q\in \ZZ\setminus\{0\}$, we
have
\[
|q\alpha-p|>\frac{c}{|q|^d}.
\]
Otherwise, $\alpha$ is called Liouville. It is easy to see that
$E^{1,0}_1=E^{0,0}_1$ and $E^{1,1}_1=E^{0,1}_1$. As shown in
\cite{Heitsch75} and \cite{Roger}, we have
\begin{itemize}
  \item $E^{0,0}_1={\mathbb R}$;
  \item $E^{0,1}_1={\mathbb R}$ if $\alpha$ is diophantine and
  $E^{0,1}_1$ is infinite dimensional if $\alpha$ is Liouville.
\end{itemize}
So when $\alpha$ is a Liouville's number, $\bar0^{0,1}_1=
\bar0^{1,1}_1\not = 0$. As a direct consequence of this fact and
\cite[Theorem D]{seq}, we obtain that, when $\alpha$ is a
Liouville's number, there exists a bundle-like metric on
$\mathbb{T}^2$ such that the associated map
$h\mapsto\cH_\infty^1(h)$, which is continuous on $[0,\infty)$ by
Corollary~\ref{c:seq} and $C^\infty$ on $(0,\infty)$, is not
$C^\infty$ at $h=0$.

\subsection{Riemannian Heisenberg manifolds} \label{s:heis}
In this Section, we discuss similar problems for adiabatic limits
associated with the Riemannian Heisenberg ma\-ni\-fold $(\Gamma
\backslash H,g)$ and the one-dimensional foliation $\mathcal F$
introduced in Section~\ref{s:Heis0} (see \cite{Heisenberg}). We
assume that $g$ corresponds to the identity matrix, and the
one-di\-men\-si\-o\-nal foliation $\cF$ is defined by the vector
field $X(1,0,0)$.

The corresponding Riemannian metric $g_h$ on $\Gamma \backslash H$
defined by (\ref{e:gh}) is given by the matrix
\[
 \begin{pmatrix}
  1 &0&0\\
  0& h^{-2}&0\\
  0&0&h^{-2}
 \end{pmatrix}, \quad h>0.
\]

By Theorem~\ref{thm:0spectrum}, it follows that the spectrum of
the Laplacian $\Delta_h $ on $0$- and $3$- forms on $M$ (with
multiplicities) is described as
\[
{\rm spec}\,\Delta^0_h ={\rm spec}\,\Delta^3_h =\Sigma_{1,h}\cup
\Sigma_{2,h},
\]
where
\begin{multline*}
\begin{aligned}
\Sigma_{1,h} & = \{\lambda_h(a,b)=4\pi^2 (a^2+h^2b^2)
: a,b\in \mathbb Z\}, \\
\Sigma_{2,h} & = \{\mu_h(c,k)=2\pi c (2k+1)h + 4\pi^2 c^{2} h^2\
\text{with mult.}\ 2c,
\end{aligned}
\\
c\in \mathbb Z^{+}, k\in \mathbb Z^{+}\cup\{0\}\}.
\end{multline*}
First, note that, for any $a\in \mathbb Z\setminus\{0\}$ and $b\in
\mathbb Z\setminus\{0\}$,
\[
\lambda_h(a,b)>4\pi^2h^2, \quad h>0,
\]
and for any $c\in \mathbb Z^{+}$ and $k\in \mathbb
Z^{+}\cup\{0\}\}$
\[
\mu_h(c,k)>4\pi^2 h^2, \quad h>0.
\]
Therefore, for any $h>0$,
\[
\lambda^0_0(h)=\lambda^3_0(h)=0, \quad
\lambda^0_1(h)=\lambda^3_1(h)>4\pi^2 h^2.
\]
Next, we see that, for any $b\in \mathbb Z\setminus\{0\}$,
\[
\lambda_h(0,b)=4\pi^2b^2h^2=O(h^2), \quad h\to0.
\]
Since we have infinitely many different branches of eigenvalues
$\lambda_h$ with $\lambda_h=O(h^2)$ as $h\to0$, we conclude that,
for any $i>0$,
\[
\lambda^0_i(h)=\lambda^3_i(h)=O(h^2), \quad h\to0.
\]

By Theorem~\ref{thm:1spectrum}, the spectrum of the Laplace
operator $\Delta_h $ on one and two forms on $M$ (with
multiplicities) has the form
\[
{\rm spec}\,\Delta^1_h ={\rm spec}\,\Delta^2_h =\Sigma_{1,h}\cup
\Sigma_{2,h}\cup \Sigma_{3,h},
\]
where
\begin{align*}
\Sigma_{1,h}  = & \{\lambda_{h,\pm}(a,b) =4\pi^{2}\left(a^2+h^2
b^2\right) + \frac{1\pm \sqrt{1+ 16\pi^{2}(a^2+h^2b^2)}}{2}\\ &
\text{with mult.}\ 2:  a,b\in \mathbb Z\},\\ \Sigma_{2,h} = &
\{\mu_h(c,k)=4\pi^2 c^{2}h^2+2\pi c (2k+1)h
\\
& \text{with mult.}\ 2c: c\in \mathbb Z^{+}, k\in \mathbb
Z^{+}\cup\{0\}\}, \\ \Sigma_{3,h}  = & \{\mu_{h,\pm}(c,k)=
4\pi^2c^2h^2+2\pi c (2k+1)h
 + \frac{1 \pm \sqrt{(4\pi c h+1)^2 +16k\pi c h}}{2}
\\ & \text{with mult.}\ 2c: c\in \mathbb Z^{+}, k\in \mathbb
Z^{+}\cup\{0\}\}.
\end{align*}

Observe that, for any $b\in \mathbb Z\setminus\{0\}$,
\[
\lambda_{h,-}(0,b) = 4\pi^{2}b^2 h^2 + \frac{1- \sqrt{1+
16\pi^{2}b^2h^2}}{2}=16\pi^{4}b^4h^4 + O(h^4), \quad h\to0,
\]
and, for any $\lambda\in {\rm spec}\,\Delta_h\setminus\{0\}$,
\[
\lambda > Ch^4, \quad h>0,
\]
with some constant $C>0$. Therefore, we have
\[
\lambda^1_0(h)=\lambda^2_0(h)=\lambda^1_1(h)=\lambda^2_1(h)=0,\quad
h>0,
\]
and, by the above argument, we obtain, for any $i>1$
\[
\lambda^1_i(h)=\lambda^2_i(h)=O(h^4), \quad h\to0.
\]

We now turn to the differentiable spectral sequence. By a
straightforward computation, one can show that $E_3=E_\infty$, all
the terms $\widehat{E}^r_1$ are infinite-dimensional and, for the
basic cohomology, we have
\[
E^{0,0}_2={\mathbb R}, \quad E^{1,0}_2={\mathbb R}, \quad
E^{2,0}_2=C^\infty(S^1).
\]
So we get that, in this case, for $r=0$ and $r=3$,
\begin{align*}
\dim{\widehat E}_1^r&= \sharp\,\left\{i\ \left|\ \lambda_i^r(h)=
O\left(h^2\right) \quad\text{as}\quad h\to 0\right.\right\}=\infty\;,\\
\dim E_k^r&= \sharp\,\left\{i\ \left|\ \lambda_i^r(h)=
O\left(h^{2k}\right) \quad\text{as}\quad h\to
0\right.\right\}=1\;,\quad k\geq2\;.
\end{align*}
and, for $r=1$ and $r=2$
\begin{align*}
\dim{\widehat E}_1^r&= \sharp\,\left\{i\ \left|\ \lambda_i^r(h)=
O\left(h^2\right) \quad\text{as}\quad h\to 0\right.\right\}=\infty\;,\\
\dim E_2^r&= \sharp\,\left\{i\ \left|\ \lambda_i^r(h)=
O\left(h^{4}\right) \quad\text{as}\quad h\to
0\right.\right\}=\infty\;,\\
\dim E_k^r&= \sharp\,\left\{i\ \left|\ \lambda_i^r(h)=
O\left(h^{2k}\right) \quad\text{as}\quad h\to
0\right.\right\}=2\;,\quad k\geq3\;.
\end{align*}
A more precise information can be obtained from the consideration
of the corresponding eigenspaces that will be discussed elsewhere.

\begin{rem}
It is quite possible that both the asymptotic formula of
Theorem~\ref{t:trace} and the investigation of small eigenvalues
given in this Section can be extended to the differential form
Laplace operator on an arbitrary Riemannian Heisenberg manifold.
Nevertheless, we believe that many essentially new features of
adiabatic limits on Riemannian Heisenberg manifolds can be already
seen in the particular cases, which were considered in this paper,
and we don't expect anything rather different in the general case.
\end{rem}

\end{document}